\documentclass[12pt]{article}
\usepackage{latexsym, amssymb}
\usepackage{dsfont}

\DeclareMathAlphabet\gothic{U}{euf}{m}{n}

\makeatletter
\def\eqnarray{\stepcounter{equation}\let\@currentlabel=\theequation
\global\@eqnswtrue
\tabskip\@centering\let\\=\@eqncr
$$\halign to \displaywidth\bgroup\hfil\global\@eqcnt\z@
  $\displaystyle\tabskip\z@{##}$&\global\@eqcnt\@ne
  \hfil$\displaystyle{{}##{}}$\hfil
  &\global\@eqcnt\tw@ $\displaystyle{##}$\hfil
  \tabskip\@centering&\llap{##}\tabskip\z@\cr}

\def\endeqnarray{\@@eqncr\egroup
      \global\advance\c@equation\m@ne$$\global\@ignoretrue}

\def\@yeqncr{\@ifnextchar [{\@xeqncr}{\@xeqncr[5pt]}}
\makeatother

\begin{document}
\bibliographystyle{tom}

\newtheorem{lemma}{Lemma}[section]
\newtheorem{thm}[lemma]{Theorem}
\newtheorem{cor}[lemma]{Corollary}
\newtheorem{voorb}[lemma]{Example}
\newtheorem{rem}[lemma]{Remark}
\newtheorem{prop}[lemma]{Proposition}
\newtheorem{stat}[lemma]{{\hspace{-5pt}}}

\newenvironment{remarkn}{\begin{rem} \rm}{\end{rem}}
\newenvironment{exam}{\begin{voorb} \rm}{\end{voorb}}

\newcommand{\gota}{\gothic{a}}
\newcommand{\gotb}{\gothic{b}}
\newcommand{\gotc}{\gothic{c}}
\newcommand{\gote}{\gothic{e}}
\newcommand{\gotf}{\gothic{f}}
\newcommand{\gotg}{\gothic{g}}
\newcommand{\gothh}{\gothic{h}}
\newcommand{\gotk}{\gothic{k}}
\newcommand{\gotm}{\gothic{m}}
\newcommand{\gotn}{\gothic{n}}
\newcommand{\gotp}{\gothic{p}}
\newcommand{\gotq}{\gothic{q}}
\newcommand{\gotr}{\gothic{r}}
\newcommand{\gots}{\gothic{s}}
\newcommand{\gotu}{\gothic{u}}
\newcommand{\gotv}{\gothic{v}}
\newcommand{\gotw}{\gothic{w}}
\newcommand{\gotz}{\gothic{z}}
\newcommand{\gotA}{\gothic{A}}
\newcommand{\gotB}{\gothic{B}}
\newcommand{\gotG}{\gothic{G}}
\newcommand{\gotL}{\gothic{L}}
\newcommand{\gotS}{\gothic{S}}
\newcommand{\gotT}{\gothic{T}}

\newcounter{teller}
\renewcommand{\theteller}{\Roman{teller}}
\newenvironment{tabel}{\begin{list}%
{\rm \bf \Roman{teller}.\hfill}{\usecounter{teller} \leftmargin=1.1cm
\labelwidth=1.1cm \labelsep=0cm \parsep=0cm}
                      }{\end{list}}

\newcounter{tellerr}
\renewcommand{\thetellerr}{(\roman{tellerr})}
\newenvironment{subtabel}{\begin{list}%
{\rm  (\roman{tellerr})\hfill}{\usecounter{tellerr} \leftmargin=1.1cm
\labelwidth=1.1cm \labelsep=0cm \parsep=0cm}
                         }{\end{list}}

\newcommand{\Ni}{\mathds{N}}
\newcommand{\Qi}{\mathds{Q}}
\newcommand{\Ri}{\mathds{R}}
\newcommand{\Ci}{\mathds{C}}
\newcommand{\N}{\mathds{N}}
\newcommand{\Q}{\mathds{Q}}
\newcommand{\R}{\mathds{R}}
\newcommand{\C}{\mathds{C}}
\newcommand{\Ti}{\mathds{T}}
\newcommand{\Zi}{\mathds{Z}}
\newcommand{\Fi}{\mathds{F}}

\newcommand{\proof}{\mbox{\bf Proof} \hspace{5pt}} 
\newcommand{\remark}{\mbox{\bf Remark} \hspace{5pt}}
\newcommand{\vertspace}{\vskip10.0pt plus 4.0pt minus 6.0pt}

\newcommand{\simh}{{\stackrel{{\rm cap}}{\sim}}}
\newcommand{\ad}{{\mathop{\rm ad}}}
\newcommand{\Ad}{{\mathop{\rm Ad}}}
\newcommand{\Aut}{\mathop{\rm Aut}}
\newcommand{\arccot}{\mathop{\rm arccot}}
\newcommand{\capp}{{\mathop{\rm cap}}}
\newcommand{\rcapp}{{\mathop{\rm rcap}}}
\newcommand{\diam}{\mathop{\rm diam}}
\newcommand{\divv}{\mathop{\rm div}}
\newcommand{\codim}{\mathop{\rm codim}}
\newcommand{\RRe}{\mathop{\rm Re}}
\newcommand{\IIm}{\mathop{\rm Im}}
\newcommand{\Tr}{{\mathop{\rm Tr}}}
\newcommand{\Vol}{{\mathop{\rm Vol}}}
\newcommand{\card}{{\mathop{\rm card}}}
\newcommand{\supp}{\mathop{\rm supp}}
\newcommand{\sgn}{\mathop{\rm sgn}}
\newcommand{\essinf}{\mathop{\rm ess\,inf}}
\newcommand{\esssup}{\mathop{\rm ess\,sup}}
\newcommand{\Int}{\mathop{\rm Int}}
\newcommand{\Leibniz}{\mathop{\rm Leibniz}}
\newcommand{\lcm}{\mathop{\rm lcm}}
\newcommand{\loc}{{\rm loc}}
\newcommand{\ran}{ {\mathcal{R}}}
\newcommand{\nul}{ {\mathcal{N}}}
\newcommand{\Hol}{\mathop{\rm Hol}}
\newcommand{\proj}{\mathop{\rm proj}}
\newcommand{\Lipp}{\mathop{\rm Lip}}

\newcommand{\mod}{\mathop{\rm mod}}
\newcommand{\spann}{\mathop{\rm span}}
\newcommand{\one}{\mathds{1}}

\hyphenation{groups}
\hyphenation{unitary}

\newcommand{\tfrac}[2]{{\textstyle \frac{#1}{#2}}}

\newcommand{\cb}{{\cal B}}
\newcommand{\cc}{{\cal C}}
\newcommand{\cd}{{\cal D}}
\newcommand{\ce}{{\cal E}}
\newcommand{\cf}{{\cal F}}
\newcommand{\ch}{{\cal H}}
\newcommand{\ci}{{\cal I}}
\newcommand{\ck}{{\cal K}}
\newcommand{\cl}{{\cal L}}
\newcommand{\cm}{{\cal M}}
\newcommand{\co}{{\cal O}}
\newcommand{\cs}{{\cal S}}
\newcommand{\ct}{{\cal T}}
\newcommand{\cw}{{\cal W}}
\newcommand{\cx}{{\cal X}}
\newcommand{\cy}{{\cal Y}}
\newcommand{\cz}{{\cal Z}}

\newlength{\hightcharacter}
\newlength{\widthcharacter}
\newcommand{\covsup}[1]{\settowidth{\widthcharacter}{$#1$}\addtolength{\widthcharacter}{-0.15em}\settoheight{\hightcharacter}{$#1$}\addtolength{\hightcharacter}{0.1ex}#1\raisebox{\hightcharacter}[0pt][0pt]{\makebox[0pt]{\hspace{-\widthcharacter}$\scriptstyle\circ$}}}
\newcommand{\cov}[1]{\settowidth{\widthcharacter}{$#1$}\addtolength{\widthcharacter}{-0.15em}\settoheight{\hightcharacter}{$#1$}\addtolength{\hightcharacter}{0.1ex}#1\raisebox{\hightcharacter}{\makebox[0pt]{\hspace{-\widthcharacter}$\scriptstyle\circ$}}}
\newcommand{\scov}[1]{\settowidth{\widthcharacter}{$#1$}\addtolength{\widthcharacter}{-0.15em}\settoheight{\hightcharacter}{$#1$}\addtolength{\hightcharacter}{0.1ex}#1\raisebox{0.7\hightcharacter}{\makebox[0pt]{\hspace{-\widthcharacter}$\scriptstyle\circ$}}}

\thispagestyle{empty}

\vspace*{1cm}
\begin{center}
{\large\bf Square roots of perturbed subelliptic operators \\[10pt]
on Lie groups} \\[5mm]
\large  Lashi Bandara, 
A.F.M. ter Elst and Alan McIntosh

\end{center}

\vspace{5mm}

\begin{list}{}{\leftmargin=1.8cm \rightmargin=1.8cm \listparindent=10mm 
   \parsep=0pt}
\item
\small
{\sc Abstract}.
We solve the Kato square root problem for bounded measurable 
perturbations of subelliptic operators on connected Lie groups.
The subelliptic operators are divergence form operators with 
complex bounded coefficients, which may have lower order
terms.
In this general setting we deduce inhomogeneous estimates.
In case the group is nilpotent and the subelliptic operator is 
pure second order, then we prove stronger homogeneous estimates.
Furthermore, we prove Lipschitz stability of the estimates under 
small perturbations of the coefficients.

\end{list}

\let\thefootnote\relax\footnotetext{
\begin{tabular}{@{}l}
\end{tabular}}

\section{Introduction} \label{Ssubkato1}

The Kato problem in $\Ri^d$ was a long standing problem which was 
 solved by Auscher, Hofmann, Lacey, McIntosh and Tchamitchian
\cite{AHLMT} in 2002.
The papers of Hofmann \cite{Hof1} and McIntosh \cite{McI4} and
the book by Auscher and Tchamitchian \cite{AT2} 
 provide a  narrative of the resolution of Kato's
conjecture. 
This problem was recast in terms of the functional calculus of a first-order system
by Axelsson, Keith and McIntosh in \cite{AKM} and \cite{AKM2}, which together provide
a unified first-order framework for recovering and extending some results concerning the harmonic analysis 
of strongly elliptic operators.
A version of the Kato problem was presented by Morris for manifolds with 
exponential growth in \cite{Morris1}, and another version on metric measure spaces 
with the doubling property was
presented by Bandara in \cite{Band1}.
The main aim of this paper is to present a solution to the Kato problem
for subelliptic operators on Lie groups.

Let $G$ be a connected Lie group with Lie algebra $\gotg$
and (left) Haar measure $\mu$.
All integration on $G$ is with respect to the (left) Haar measure
and the norm is the $L_2$-norm,
unless stated otherwise.
Let $a_1,\ldots,a_m$ be an algebraic basis for $\gotg$, that is, 
an independent set which generates $\gotg$.
Let $L$ be the left regular representation in $L_2(G)$.
So $(L(x) f)(y) = f(x^{-1} y)$ for all $x \in G$, $f \in L_2(G)$
and a.e.\ $y \in G$.
For all $k \in \{ 1,\ldots,m \} $ let $A_k$ be the 
infinitesimal generator of the one-parameter unitary group
$t \mapsto L(\exp t a_k)$. 
Then $A_k$ is skew-adjoint.
Define the Sobolev space 
$W_{1,2}'(G) = \bigcap_{k=1}^m D(A_k)$ with norm such that 
\[
\|f\|_{W_{1,2}'(G)}^2
= \|f\|^2 + \sum_{k=1}^m \|A_k f\|^2
.  \]
Then $W_{1,2}'(G)$ is a Hilbert space since $A_k$ is closed for all~$k$.
Note that $W_{1,2}'(G)$ depends on the choice of the algebraic basis
and if confusion is possible, then we write $W_{1,2}'(G,a)$.
Clearly $C_c^\infty(G) \subset W_{1,2}'(G)$, so $W_{1,2}'(G)$ is dense in $L_2(G)$.

Next, for all $k,l \in \{ 1,\ldots,m \} $ let $b_{kl}, b_k, b_k', b_0 \in L_\infty(G)$.
Assume there exists a constant $\kappa > 0$ such that the following G\aa rding inequality holds:  
\begin{eqnarray}
\lefteqn{
\RRe \Big( \sum_{k,l=1}^m (b_{kl} \, A_l u, A_k u) 
   +\sum_{l=1}^m (b_{l} \, A_l u, u)
   +\sum_{k=1}^m (b_{k}' \, u, A_k u)
   +(b_0 \, u,u) \Big)
} \hspace*{100mm} \label{eSsubkato1;20}  \\
& \geq & \kappa \Big( \sum_{k=1}^m \|A_k u\|^2 + \|u\|^2 \Big)
\nonumber 
\end{eqnarray}
for all $u \in W_{1,2}'(G)$.
Then 
the inhomogeneous divergence form operator
\[
H_I =
- \sum_{k,l=1}^m A_k \, b_{kl} \, A_l 
   + \sum_{l=1}^m b_l \, A_l
   - \sum_{k=1}^m A_k \, b_k'
   + b_0 \, I
\]
is a maximal accretive operator on $L_2(G)$ of type $\omega$ for some 
$\omega \in [0,\frac{\pi}{2})$. So  $- H_I$
generates a bounded semigroup on $L_2(G)$.

It is easy to see that (\ref{eSsubkato1;20}) holds under the homogeneous G\aa rding inequality
\begin{equation}
\RRe \sum_{k,l=1}^m (b_{kl} \, A_l u, A_k u) \geq \kappa \sum_{k=1}^m \|A_k u\|^2
\label{homgarding}
\end{equation}
with a possible change in $\kappa$, provided a large enough positive constant is added to $b_0$.
For example, if 
\[
\RRe \sum_{k,l=1}^m b_{kl} \, \xi_k \, \overline{\xi_l} 
\geq \kappa \, |\xi|^2
\]
a.e.\ for all $\xi\in \Ci^m$, then (\ref{homgarding}) is valid.

Furthermore, let $b \in L_\infty(G)$ and suppose there exists a constant 
$\kappa_1 > 0$ such that $\RRe b \geq \kappa_1$ a.e.
Then $b H_I$ is an $\tilde\omega$-sectorial operator  on  $L_2(G)$ for some $\tilde\omega < \pi$,
so it has a unique square root $\sqrt{b H_I}$ which is $\frac12\tilde\omega$-sectorial and 
satisfies $(\sqrt{b H_I})^2= bH_I$.

The main theorems of this paper are as follows.
The first one is the solution of the inhomogeneous Kato
problem for subelliptic operators.

\begin{thm} \label{tsubkato101}
Let $G$ be a connected Lie group and suppose that 
$a_1,\ldots,a_m$ is an algebraic basis for the Lie algebra $\gotg$ of $G$.
Let 
\[
H_I =
- \sum_{k,l=1}^m A_k \, b_{kl} \, A_l 
   + \sum_{l=1}^m b_l \, A_l
   - \sum_{k=1}^m A_k \, b_k'
   + b_0 \, I
\]
be a divergence form operator with bounded measurable coefficients
satisfying the ellipticity condition {\rm (\ref{eSsubkato1;20})}.
Let $b \in L_\infty(G)$ and suppose there exists a constant $\kappa_1 > 0$ such that 
$\RRe b \geq \kappa_1$ a.e.
Then $D(\sqrt{b H_I}) = W_{1,2}'(G)$ and there exist $c,C > 0$ such that 
\[
c \, (\|\sqrt{b H_I} u\| + \|u\|)
\leq \|u\| + \sum_{k=1}^m \|A_k u\|
\leq C \, (\|\sqrt{b H_I} u\| + \|u\|)
\]
for all $u \in W_{1,2}'(G)$.
\end{thm}

For connected nilpotent Lie groups, or more generally, for Lie groups $G$
which are the local direct product of a connected compact Lie group $K$
and a connected nilpotent Lie group $N$, a homogeneous result is also valid.
(To say that $G$ is 
the local direct product of $K$
and  $N$ means that
$G = K \cdot N$ and $K \cap N$ is discrete.
Equivalently, the Lie algebra of $G$ is the direct product of the 
Lie algebra of $K$ and the Lie algebra of $N$.) 

\begin{thm} \label{tsubkato102}
Let $G$ be the local direct product of a connected compact Lie group 
and a connected nilpotent Lie group.
Let 
\[
H = - \sum_{k,l=1}^m A_k \, b_{kl} \, A_l 
\]
be a homogeneous divergence form operator with bounded measurable coefficients
satisfying the subellipticity condition {\rm (\ref{homgarding})}.
Let $b \in L_\infty(G)$ and suppose there exists a constant $\kappa_1 > 0$ such that 
$\RRe b \geq \kappa_1$ a.e.
Then $D(\sqrt{b H}) = W_{1,2}'(G)$ and there exist $c,C > 0$ such that 
\[
c \, \|\sqrt{b H} u\| 
\leq \sum_{k=1}^m \|A_k u\|
\leq C \, \|\sqrt{b H} u\| 
\]
for all $u \in W_{1,2}'(G)$.
\end{thm}
 
Lie groups of the above mentioned type necessarily satisfy the doubling property. We do not know whether the doubling property itself  implies homogeneous bounds, except in
the special case when the algebraic basis $a_1,\ldots,a_m$ is actually a 
vector space basis.  In other words, when $H$ is strongly elliptic rather 
than subelliptic, then the conclusion of Theorem~\ref{tsubkato102} does
hold on all Lie groups with polynomial growth. 
 For vector space bases we drop the prime and write $W_{1,2}(G) = W_{1,2}'(G)$
and $W_{1,2}(G,a) = W_{1,2}'(G,a)$. 

\begin{thm} \label{tsubkato303}
Let $G$ be a connected Lie group with polynomial growth and suppose that 
$a_1,\ldots,a_m$ is a vector space basis for the Lie algebra $\gotg$ of $G$.
Let 
\[
H = - \sum_{k,l=1}^m A_k \, b_{kl} \, A_l 
\]
be a strongly elliptic
homogeneous divergence form operator with bounded measurable coefficients
satisfying the ellipticity condition {\rm (\ref{homgarding})}.
Let $b \in L_\infty(G)$ and suppose there exists a constant $\kappa_1 > 0$ such that 
$\RRe b \geq \kappa_1$ a.e.
Then $D(\sqrt{b H}) = W_{1,2}(G)$ and there exist $c,C > 0$ such that 
\[
c \, \|\sqrt{b H} u\| 
\leq \sum_{k=1}^m \|A_k u\|
\leq C \, \|\sqrt{b H} u\| 
\]
for all $u \in W_{1,2}(G)$.
\end{thm}
 
In Section~\ref{Ssubkato3} we prove the homogeneous bounds, first those in Theorem \ref{tsubkato102}. 
The algebraic basis provides a canonical distance $d$ that is
well suited for the study of subelliptic operators.
Then $(G,d,\mu)$ is a  metric measure space.  The proof is achieved by building upon the results
of Bandara \cite{Band1} who adapted the earlier framework of \cite{AKM} to the situation of metric measure spaces.

The situation in  Theorem~\ref{tsubkato303} requires more substantial innovations in the proof, involving  the structure theory of Lie groups as developed by Dungey--ter Elst--Robinson \cite{DER4}.  Note that  we cannot  directly apply 
 homogeneous estimates for second-order derivatives of the form
$\|A_k \, A_l u\| \leq C \|\Delta u\|$
where $\Delta$ is the (sub-)Laplacian, as typically used in proofs of the Kato estimates. 
Indeed these
homogeneous estimates $\|A_k \, A_l u\| \leq C \|\Delta u\|$
for second-order derivatives are {\em false} for $G$, if 
$G$ is a connected Lie group with polynomial growth
which is not a local direct product of a connected compact Lie group 
and a connected nilpotent Lie group. (See \cite{ERS2} Theorem~1.1.)

In Section 4 we turn to the proof of the inhomogeneous estimates as stated in Theorem \ref{tsubkato101}. The Haar measure is at most exponential 
in volume growth of balls, so we are in a position where we can adapt the results of Morris \cite{Morris1} who obtained Kato estimates on complete Riemannian manifolds which satisfy a similar exponential  growth condition on the volume of balls.
Indeed his methods work for an arbitrary Borel-regular measure.
 
In Section \ref{Ssubkato5} we consider some variants of the inhomogeneous results, while in Section \ref{Ssubkato6} we state and prove a Lipschitz estimate of the form (in the homogeneous case)
\[
\|\sqrt{(b + \tilde b) \widetilde H} u - \sqrt{b H} u\| 
\leq C (\|\tilde b\|_\infty + \sum_{j,k}\|\tilde b_{jk}\|_\infty) \|\nabla u\|\ ,
\]
where
\[
H = - \sum_{k,l=1}^m A_k \, b_{kl} \, A_l \qquad
\mbox{and} \qquad
\widetilde H = - \sum_{k,l=1}^m A_k \, (b_{kl} + \tilde b_{kl}) \, A_l \ ,
  \]
 under small bounded perturbations $\tilde b, \tilde b_{kl}$ of the coefficients.

\section{Preliminaries} \label{Ssubkato2}

In this section we gather some background material on Lie groups and 
operator theory that will be used throughout the paper.

\subsection{Lie groups}  \label{Ssubkato2.1}

We use the notation as in the introduction.
In particular, $a_1,\ldots,a_m$ is an algebraic basis for the 
Lie algebra $\gotg$ of a connected Lie group $G$.

For all $k \in \{ 1,\ldots,m \} $ and $x \in G$ define
$X_k|_x \in T_xG$ by 
\[
X_k|_x f
= \frac{d}{dt} f((\exp t a_k) x) \Big|_{t=0}
.  \]
Then $X_k$ is a smooth right invariant vector field on $G$.
Note that $A_k f = - X_k f$ for all $f \in C_c^\infty(G)$.

The space $L_2(G,\Ci^m)$ has a natural inner product.
Define the unbounded operator $\nabla$ from $L_2(G)$ into 
$L_2(G,\Ci^m)$ by $D(\nabla) = W_{1,2}'(G)$ and 
\[
(\nabla u)(x)
= ((A_1 u)(x), \ldots, (A_m u)(x))
\]
for a.e.\ $x \in G$.
Then $\nabla$ is densely defined and closed, since $A_k$ is closed
for all $k \in \{ 1,\ldots,m \} $.
We denote its adjoint by $-\divv$. 
Thus $\divv = - \nabla^*$.

Define the sesquilinear form $J \colon W_{1,2}'(G) \times W_{1,2}'(G) \to \Ci$ by
\[
J[f,g] = (\nabla f, \nabla g)
.  \]
Then $J$ is a closed positive symmetric form.
Let $\Delta$ be the self-adjoint operator associated with $J$.
We call $\Delta$ the {\bf sub-Laplacian} (associated with the algebraic basis
$a_1,\ldots,a_m$).
Clearly $\Delta = - \divv \, \nabla$.
There is another simple identity for $\Delta$.

\begin{prop} \label{psubkato202}
One has $\Delta = - \sum_{k=1}^m A_k^2$ with
$D(\Delta) = \bigcap_{k=1}^m D(A_k^2)$.
Moreover, $D(\Delta) = \bigcap_{k,l=1}^m D(A_k \, A_l)$
and there exists a constant $C_1 > 0$ such that 
\begin{equation}
\sum_{k,l=1}^m \|A_k \, A_l u\|^2 \leq C_1 (\|\Delta u\|^2 + \|u\|^2)
\label{epsubkato202;1}
\end{equation}
for all $u \in D(\Delta)$.

Finally, if $G$ is the local direct product of a connected compact Lie group 
and a connected nilpotent Lie group, then there exists a constant $C_2 > 0$
such that
\[
\frac{1}{m} \, \|\Delta u\|^2
\leq \sum_{k,l=1}^m \|A_k \, A_l u\|^2 
\leq C_2 \, \|\Delta u\|^2
\] 
for all $u \in D(\Delta)$. 
\end{prop}
\proof\
It is proved in ter Elst--Robinson \cite{ER1} Theorem~7.2.I
that $\Delta_0 := \Delta|_{\cap_{k,l=1}^m D(A_k A_l)}$ is a 
closed operator and in the proof it is shown that $\Delta_0 = \Delta$
(with same domains).
Moreover, ter Elst--Robinson \cite{ER8} Theorem~3.3.III
gives $\bigcap_{k,l=1}^m D(A_k \, A_l) = \bigcap_{k=1}^m D(A_k^2)$.
Then estimate (\ref{epsubkato202;1}) follows from the closed graph theorem.
The final equivalence of seminorms is proved in 
ter Elst--Robinson--Sikora \cite{ERS2} Proposition~4.1.\hfill$\Box$

\vertspace

Define 
\[
W_{1,2}'(G,\Ci^m) 
= \{ x \mapsto (u_1(x),\ldots,u_m(x)) : u_1,\ldots,u_m \in W_{1,2}'(G) \}
\subset L_2(G,\Ci^m)
.  \]
Next, the space $L_2(G,\Ci^{m^2})$ has a natural inner product.
Define the unbounded operator $\widetilde \nabla$ from $L_2(G,\Ci^m)$ into 
$L_2(G,\Ci^{m^2})$ by $D(\widetilde \nabla) = W_{1,2}'(G,\Ci^m)$ and 
\begin{equation}\label{tildenabla}
(\widetilde \nabla f)(x)
= \Big( (A_k u_l)(x) \Big)_{kl}
\end{equation}
for a.e.\ $x \in G$, if 
$f(x) = (u_1(x),\ldots,u_m(x))$ for a.e.\ $x \in G$.

\begin{remarkn}
To satisfy our readers with a 
geometric appetite, we make the following geometric remark.
Recall the notion of a connection $\nabla$ 
over a vector bundle $V$. This is an operator
$\nabla \colon C_{\infty}(V) \to C_{\infty}(T^\ast G \otimes V)$ 
satisfying $\nabla_{fX}(Y) = f\nabla_{X} Y$
and $\nabla_{X}(fY) = X(f) Y  + f \nabla_{X} Y$
for all $f \in C_{\infty}(G)$,
$X \in C_{\infty}(T^\ast G)$ and
$Y \in C_{\infty}(V)$.
Now, given a sub-bundle
$E \subset T^\ast G$, we can define a
{\bf sub-connection}
$\nabla \colon C_{\infty}(V) \to C_{\infty}(E \otimes V)$
to mean that it satisfies
the above properties but with the
condition on $X$ being that 
$X \in C_{\infty}(E)$.
Our philosophy in this paper
stems from an observation that
we can construct a sub-bundle
$E$ on which our subelliptic
operators are strongly elliptic.
Since Lie groups are parallelisable,
and for the benefit of a wider
audience, we refrain from using
the language of vector bundles
in this paper. 
However, we refer
the reader to \cite{BMC},
where the first and third authors
provide a solution to a version
of the Kato square root problem
for certain uniformly
elliptic second order operators
over a class of vector bundles.
\end{remarkn}

We also need a subelliptic distance on~$G$.
Let $\gamma \colon [0,1] \to G$ be an absolutely continuous path such that 
$\dot\gamma(t) \in \spann \{ X_1|_{\gamma(t)}, \ldots, X_m|_{\gamma(t)} \}  $ 
for a.e.\ $t \in [0,1]$.
Define the {\bf length} of $\gamma$ by 
\[
\ell(\gamma) = \int_0^1  \Big( \sum_{k=1}^m |\gamma^k(t)|^2 \Big)^{1/2} \, dt
\in [0,\infty]
\]
if $\dot\gamma(t) = \sum_{k=1}^m \gamma^k(t) \, X_k|_{\gamma(t)}$
for a.e.\ $t \in [0,1]$.
Since $G$ is connected it follows from a theorem of Carath\'eodory
\cite{Carath} that for all $x,y \in G$ there exists such a path 
$\gamma$ with finite length and $\gamma(0) = x$ and $\gamma(1) = y$.
If $x,y \in G$ then we define the {\bf distance} $d(x,y)$ between $x$ and $y$ to be the 
infimum of the length of all such paths with $\gamma(0) = x$ and $\gamma(1) = y$.
Then $d$ is a metric on $G$.
For all $x \in G$ and $r > 0$, let $B(x,r) = \{ y \in G : d(x,y) < r \} $.
If $x = e$, the identity element of $G$, then we write $B(r) = B(e,r)$.
Then $B(x,r) = B(r) x$ for all $x \in G$.

\begin{prop} \label{psubkato203}
\mbox{}
\begin{tabel}
\item \label{psubkato203-0.5}
The topology on $G$ is the same as the 
topology associated with $d$.
In particular, the open balls are measurable.
\item \label{psubkato203-1}
The metric space $(G,d)$ is complete and the closed balls $\overline{B(x,r)}$ are compact.
\item \label{psubkato203-2}
There exist $c,C > 0$ and $D' \in \Ni$ such that 
$c \, r^{D'} \leq \mu(B(r)) \leq C \, r^{D'}$
for all $r \in (0,1]$.
\item \label{psubkato203-3}
There exist $C > 0$ and $\lambda \geq 0$ such that 
$\mu(B(r)) \leq C \, e^{\lambda r}$ for all $r \geq 1$.
\end{tabel}
\end{prop}
\proof\
For a proof of Statement~\ref{psubkato203-0.5}, see \cite{VSC} Proposition~III.4.1.

Statement~\ref{psubkato203-1} follows from the discussion in Section~III.4 
in \cite{VSC} and the fact that every locally compact metric space is complete.

Statement~\ref{psubkato203-2} is a consequence of 
Nagel--Stein--Wainger \cite{NSW} Theorems~1 and~4.

The last statement is proved in Guivarc'h \cite{Guiv} Th\'eor\`eme~II.3.\hfill$\Box$

\vertspace

A {\bf metric measure space} $(\cx,\rho,\nu)$ is a set $\cx$ with a metric $\rho$
and a measure $\nu$ on the Borel $\sigma$-algebra of $\cx$ induced by the 
metric on $\cx$.
The metric measure space is said to have the {\bf doubling property}
if there exists a constant $c > 0$ such that 
$0 < \nu(B(x,2r)) \leq c \, \nu(B(x,r)) < \infty$ for all $x \in \cx$ and $r > 0$.
It is called {\bf locally exponentially doubling} if 
there exist $\kappa, \lambda \geq 0$ and $C \geq 1$ such that  
\[
0 < \nu(B(x,tr)) \leq C \, t^{\kappa} \, e^{\lambda tr} \nu(B(x,r)) < \infty
\]
for all $x \in \cx$, $r > 0$ and $t \geq 1$.
It follows from Proposition~\ref{psubkato203}.\ref{psubkato203-0.5}
that the triple $(G,d,\mu)$ is a metric measure space.
We say the Lie group $G$ has {\bf polynomial growth} if there exist $D \in \Ni$
and $c > 0$ 
such that $\mu(B(r)) \leq c \, r^D$ for all $r \geq 1$.
 Recall that the {\bf modular function} on $G$, 
which we denote by $\delta$ throughout this paper, is the function
$\delta \colon G \to (0,\infty)$ such that 
$\mu(U \, x) = \delta(x) \, \mu(U)$ for every Borel measurable $U \subset G$
and $x \in G$.
It is a continuous homomorphism. 
The group $G$ is called {\bf unimodular} if $\delta(x) = 1$ for all $x \in G$.

\begin{prop} \label{psubkato203.5}
\mbox{}
\begin{tabel}
\item \label{psubkato203.5-1}
The metric measure space $(G,d,\mu)$ is locally exponentially doubling.
\item \label{psubkato203.5-4}
The Lie group $G$ has polynomial growth if and only if 
the metric measure space $(G,d,\mu)$ has the doubling property.
If $G$ has polynomial growth then 
there exist $c,C > 0$ and $D \in \Ni_0$ such that 
$c \, r^D \leq \mu(B(x,r)) \leq C \, r^D$
for all $x \in G$ and $r \geq 1$.
\end{tabel}
\end{prop}
\proof\
`\ref{psubkato203.5-1}'. 
It follows from Statements~\ref{psubkato203-2} and \ref{psubkato203-3}
of Proposition~\ref{psubkato203} that there exist $C \geq 1$, $\lambda \geq 0$ and $D' \in \Ni$
such that $0 < \mu(B(tr)) \leq C \, t^{D'} \, e^{\lambda t r} \, \mu(B(r))$
for all $r > 0$ and $t \geq 1$.
Then for all $x \in G$, $r > 0$ and $t \geq 1$ one has 
\begin{eqnarray*}
0 < \mu(B(x,tr)) 
= \mu(B(tr) x)
& = & \delta(x) \, \mu(B(tr))  \\
& \leq & C \, \delta(x) \, t^{D'} \, e^{\lambda t r} \, \mu(B(r))
= C \, t^{D'} \, e^{\lambda t r} \, \mu(B(x,r))
.  
\end{eqnarray*}
So the metric measure space $(G,d,\mu)$ is locally exponentially doubling.

`\ref{psubkato203.5-4}'. 
Clearly volume doubling implies polynomial growth.
Conversely, suppose that $G$ has polynomial growth.
Then it follows from Guivarc'h \cite{Guiv} Th\'eor\`eme~II.3
that there exist $c,C > 0$ and $D \in \Ni_0$ such that 
$c \, r^D \leq \mu(B(r)) \leq C \, r^D$
for all $r \geq 1$.
Moreover, $G$ is unimodular by \cite{Guiv} Lemma~1.3.
Therefore $\mu(B(x,r)) = \mu(B(r) x) = \mu(B(r))$ for all $x \in G$ and $r > 0$
and the volume estimate follows.
Together with Proposition~\ref{psubkato203}.\ref{psubkato203-2} it follows that there
exists a constant $c' > 0$ such that 
$\mu(B(2r)) \leq c' \, \mu(B(r))$ for all $r > 0$.
Therefore the metric measure space $(G,d,\mu)$ has the doubling property.\hfill$\Box$

\vertspace

If $\eta$ is a scalar valued Lipschitz function $\eta$ on a metric space 
$(\cx,\rho)$  without isolated points, then 
we denote by $\Lipp\eta$ the pointwise Lipschitz constant
\[
(\Lipp\eta)(x) 
= \limsup_{y \to x} \frac{|\eta(x) - \eta(y)|}{\rho(x,y)}
\]
whenever $x \in \cx$.

Let $L^{(1)}$ be the left regular representation in $L_1(G)$ and 
for all $k \in \{ 1,\ldots,m \} $ let $A^{(1)}_k$ be the 
infinitesimal generator of the one-parameter group
$t \mapsto L^{(1)}(\exp t a_k)$.
Let $A^{(\infty)}_k = - (A^{(1)}_k)^*$ be the dual operator on $L_\infty(G)$.
One has the following relation with Lipschitz functions.

\begin{prop} \label{psubkato210}
\mbox{}
\begin{tabel}
\item \label{psubkato210-1}
The space $\bigcap_{k=1}^m D(A^{(\infty)}_k)$ is the space 
of all bounded Lipschitz functions on $G$.
\item \label{psubkato210-2}
If $\eta \in \bigcap_{k=1}^m D(A^{(\infty)}_k)$, then 
$\sum_{k=1}^m (A^{(\infty)}_k \eta)^2 \leq (\Lipp \eta)^2 $ a.e.
\end{tabel}
\end{prop}
\proof\
Statement~\ref{psubkato210-1} follows as in Theorem~6.12 in \cite{Hei}.

`\ref{psubkato210-2}'.
Let $\xi \in \Ri^m$ with $|\xi| = 1$.
For all $t > 0$ define $y_t = \exp(t \sum_{k=1}^m \xi_k \, a_k)$.
Then $d(y_t,e) \leq t$ and if $t$ is small enough then $y_t \neq e$.
Let $x \in G$.
Then for all $t > 0$ with $y_t \neq e$ one has 
\begin{eqnarray*}
\sup_{y \in B(x,t) \setminus \{ x \} }
    \frac{|\eta(y) - \eta(x)|}{d(y,x)} 
& \geq & \frac{|\eta(y_t \, x) - \eta(x)|}{d(y_t \, x,x)} 
\geq \frac{|\eta(y_t \, x) - \eta(x)|}{t} 
= \Big| \frac{1}{t} \Big( (L^{(\infty)}(y_t^{-1}) - I) \eta \Big)(x) \Big|
,
\end{eqnarray*}
where $L^{(\infty)}$ is the left regular representation in $L_\infty(G)$.
Now choose $t = \frac{1}{n}$ with $n \in \Ni$ and note that 
\[
\lim_{n \to \infty} n (L^{(\infty)}(y_{1/n}^{-1}) - I) \eta
= - \sum_{k=1}^m \xi_k \, A^{(\infty)}_k \eta
\]
weakly$^*$ in $L_\infty(G)$.
Therefore 
\[
|\Lipp \eta| 
\geq \Big| \sum_{k=1}^m \xi_k \, A^{(\infty)}_k \eta \Big|
\]
a.e.
This is for all $\xi \in \Ri^m$ with $|\xi| = 1$.

Let $D$ be a countable dense subset of 
$ \{ \xi \in \Ri^m : |\xi| = 1 \} $.
Then there exists a nulset $N \subset G$ such that 
\begin{equation}
|(\Lipp \eta)(x)| 
\geq \Big| \sum_{k=1}^m \xi_k \, (A^{(\infty)}_k \eta)(x) \Big|
\label{epsubkato210;1}
\end{equation}
for all $x \in G \setminus N$ and $\xi \in D$.
Hence by continuity (\ref{epsubkato210;1}) is valid for all 
$x \in G \setminus N$ and $\xi \in \Ri^m$ with $|\xi| = 1$.
Therefore 
\[
|(\Lipp \eta)(x)|^2
\geq \sum_{k=1}^m |(A^{(\infty)}_k \eta)(x)|^2
\]
for all $x \in G \setminus N$.\hfill$\Box$

\subsection{Operator theory}  \label{Ssubkato2.2}

For the convenience of the reader, we shall present some operator
theoretic material which sits at the heart of both the
homogeneous and inhomogeneous problems.

First, we recall the theory of {\bf bisectorial operators}. 
For all $\omega \in [0,\frac{\pi}{2})$ 
define the {\bf bisector} by
\[
S_\omega = \{ \zeta \in \C : |\arg \zeta| \leq \omega \mbox{ or } 
                             |\pi - \arg \zeta| \leq \omega \mbox{ or } \zeta = 0 \} 
\]
and the {\bf open bisector} by 
\[
S_\omega^o = \{ \zeta \in \C \setminus \{ 0 \} : 
     |\arg \zeta| < \omega \mbox{ or } |\pi - \arg \zeta| < \omega \}
.  \]
Let $\ch$ be a Hilbert space. 
An operator $T \colon D(T)  \to \ch$ with $D(T) \subset \ch$
is then called {\bf $\omega$-bisectorial} 
(or bisectorial with {\bf angle of sectoriality} $\omega$)
if it is closed, $\sigma(T) \subset S_\omega$, 
and for each $\omega < \mu < \frac{\pi}{2}$, there is
a constant $C_\mu > 0$ such that $|\zeta| \, \|(\zeta \, I - T)^{-1}\| \leq C_\mu$
for all $\zeta \in \C \setminus \{ 0 \} $ satisfying 
$|\arg\zeta| \geq \mu$ and $|\pi - \arg\zeta| \geq \mu$.
 
\begin{remarkn} 
When $T$ is $\omega$-bisectorial, then $T^2$ is $2\omega$-sectorial 
(meaning that 
$\sigma(T) \subset S^+_{2\omega} 
= \{ \zeta \in \C : |\arg \zeta| \leq 2\omega \mbox{ or } \zeta = 0 \} $ 
and appropriate resolvent bounds hold) and hence $T^2$ has a unique 
$\omega$-sectorial square root $\sqrt{T^2}$. It may or may not happen that  
$D(\sqrt{T^2})=D(T)$ with homogeneous 
($\|\sqrt{T^2}u\|\simeq \|Tu\|$) or inhomogeneous 
($\|\sqrt{T^2}u\|+\|u\|\simeq \|Tu\|+\|u\|$) equivalence of norms. 
The determination of such equivalences involves studying the holomorphic  
functional calculus of $T$ and proving quadratic estimates.
\end{remarkn}

Let $T$  be an $\omega$-bisectorial operator. 
Then one has the  (possibly non-orthogonal)
decomposition $\ch = \nul(T) \oplus \overline{\ran(T)}$ by 
a variation of the proof of Theorem~3.8 in 
Cowling--Doust--McIntosh--Yagi \cite{CDMY}.
We denote by $\proj_{\nul(T)}$ the projection from $\ch$ onto 
$\nul(T)$ along this decomposition.
Bisectorial operators admit a functional calculus in the following sense \cite{ADM}.
For all $\mu \in (0,\frac{\pi}{2})$ let $\Psi(S_\mu^o)$ denote the space of all
holomorphic functions $\psi \colon S_\mu^o \to \C$ 
for which there exist $\alpha,c > 0$ such that 
\[
|\psi(\zeta)| \leq c \, \frac{|\zeta|^\alpha}{1 + |\zeta|^{2\alpha}}
\]
for all $\zeta \in \Psi(S_\mu^o)$.
If $\mu > \omega$ then for all $\psi \in \Psi(S_\mu^o)$ one can define
the bounded operator
\[
\psi(T) = \frac{1}{2\pi i} \oint_\gamma \psi(\zeta) \, (\zeta \, I - T)^{-1} \, d\zeta
,  \]
where $\gamma$ is a contour in $S_\mu^o$ enveloping $S_\omega$
parametrised anti-clockwise.
The integral here is simply defined via Riemann sums and this sum
converges absolutely as a consequence of the decay of $\psi$
coupled with the resolvent bounds of $T$.
If, in addition, there exists a constant $C > 0$ such that 
$\|\psi(T)\| \leq C \, \|\psi\|_\infty$ for all $\psi \in \Psi(S_\mu^o)$,
then we say that $T$ has a {\bf bounded holomorphic $S_\mu^o$ functional calculus}.

Define $\Hol^\infty(S_\mu^o)$ to be 
the space of all bounded functions $f \colon S_\mu^o \cup \{ 0 \} \to \C$
which are holomorphic on $S_\mu^o$.
For all $f \in \Hol^\infty(S_\mu^o)$ there exists a uniformly bounded 
sequence $(\psi_n)_{n \in \Ni}$
in $\Psi(S_\mu^o)$ which converges to $f$ on $S_\mu^o$
in the compact-open topology.
If in addition $T$ has a bounded holomorphic $S_\mu^o$ functional calculus, then
$\lim_{n \to \infty} \psi_n(T)$
exists in the strong operator topology on $\cl(\ch)$ by a 
modification of the proof of the theorem in Section~5 in McIntosh \cite{McI2},
and we define $f(T) \in \cl(\ch)$ by 
$f(T)u = \lim_n \psi_n(T)u+ f(0) \, \proj_{\nul(T)}u$ for all $u \in \ch$. 
The bounded operator $f(T)$ is indeed independent of the choice of the 
sequence $(\psi_n)_{n \in \Ni}$.
We then say that $T$ has a {\bf bounded holomorphic $H^\infty(S_\mu^o)$-functional calculus}.

Define $\chi^\pm \colon \Ci \to \Ci$ by
\[
\chi^\pm(z) 
= \left\{ \begin{array}{ll}
   1 & \mbox{if } \pm \RRe z > 0 ,  \\[5pt]
   0 & \mbox{if } \pm \RRe z \leq 0 .
          \end{array} \right.
\]
Moreover, define $\sgn = \chi^+ - \chi^-$.
Then $\chi^\pm,\sgn \in \Hol^\infty(S_\mu^o)$ for all $\mu \in (0,\frac{\pi}{2})$.

Next we recall some important facts from Axelsson--Keith--McIntosh \cite{AKM} regarding quadratic estimates.
We consider a triple $(\Gamma, B_1, B_2)$ of operators in $\ch$.
First, we quote the following hypotheses from this reference.

\begin{enumerate}
\item[(H1)]
	The operator $\Gamma \colon D(\Gamma) \to \ch$
is a closed, densely defined operator from $D(\Gamma) \subset \ch$ into~$\ch$
such that $\ran(\Gamma) \subset \nul(\Gamma)$.
So $\Gamma^2 = 0$.

\item[(H2)]
	The operators $B_1$ and $B_2$ are bounded on $\ch$.
Moreover, there are $\kappa_1,\kappa_2 > 0$ such that 
\begin{eqnarray*}
\RRe (B_1 u, u) & \geq & \kappa_1 \, \|u\|^2 \quad \mbox{for all } u \in \ran(\Gamma^*) , \mbox{ and}  \\
\RRe (B_2 u, u) & \geq & \kappa_2 \, \|u\|^2 \quad \mbox{for all } u \in \ran(\Gamma) .
\end{eqnarray*}

\item[(H3)]
	The operators $B_1$ and $B_2$ satisfy
	$B_1 \, B_2 (\ran(\Gamma)) \subset \nul(\Gamma)$
	and $B_2 \, B_1(\ran(\Gamma^*)) \subset \nul(\Gamma^*)$.
\end{enumerate}

Define $\Pi = \Gamma + \Gamma^*$ and $\Pi_B = \Gamma + B_1 \, \Gamma^* \, B_2$.
Then $\Pi$ is self-adjoint.
If 
\begin{equation}
\omega = \frac{1}{2} \Big( \arccos \frac{\kappa_1}{\|B_1\|} + \arccos \frac{\kappa_2}{\|B_2\|} \Big)
\label{eSsubkato2;12}
\end{equation}
then it follows from Proposition~2.5 in \cite{AKM} that 
$\Pi_B$ is $\omega$-bisectorial.
Hence for all $t > 0$ one can define the 
operator $Q^B_t = i \Pi_B (I + t^2 \Pi_B^2)^{-1} \in \cl(\ch)$.

The following proposition highlights the connection 
between the harmonic analysis and bounded holomorphic functional calculus.

\begin{thm}[Kato square root type estimate] \label{tsubkato211}
Suppose the triple $(\Gamma,B_1,B_2)$
satisfies {\rm (H1)--(H3)} and the operator
$\Pi_B$ satisfies the quadratic estimate
\[
\int_{0}^\infty \|Q_t^B u\|^2\ \frac{dt}{t}
	\simeq \|u\|^2
\] 
for all $u \in \overline{\ran(\Pi_B)}$.
Then the following hold:
\begin{tabel}
\item \label{tsubkato211-0.4}
For all $\mu \in (\omega,\frac{\pi}{2})$ the operator $\Pi_B$ has a bounded 
holomorphic $S_\mu^o$ functional calculus, where $\omega$ is as 
in~{\rm (\ref{eSsubkato2;12})}.
\item \label{tsubkato205-1}
The Hilbert space admits the spectral decomposition
\[
\ch = \nul(\Pi_B) \oplus \overline{\ran(\Pi_B)}
= \nul(\Pi_B) \oplus \ran(\chi^+(\Pi_B)) \oplus \ran(\chi^-(\Pi_B))
,  \]
where the sums are in general not orthogonal.
\item \label{tsubkato205-2}
One has $D(\Pi_B) = D(\sqrt{\Pi_B^2})$ and
\[
\|\Gamma u\| + \|B_1 \, \Gamma^* \, B_2 \, u\| 
\simeq \|\Pi_B \, u\|
\simeq \|\sqrt{\Pi_B^2} \, u\| 
\]
for all $u \in D(\Pi_B)$.
\end{tabel}
\end{thm}
\proof\
`\ref{tsubkato211-0.4}'.
This is contained in the proof of Proposition~4.8 in 
\cite{AKM}. 

`\ref{tsubkato205-1}'. 
The first equality follows from \cite{AKM} Proposition~2.2.
The second equality follows from the bounded $H^\infty$ functional calculus
of $\Pi_B$.

`\ref{tsubkato205-2}'. 
This follows from Lemma~4.2 in \cite{AKM} and again functional calculus.
Essentially it uses the fact that $\sgn(\Pi_B)$ is bounded 
(by \ref{tsubkato211-0.4}), along with  the identities 
$\sqrt{\Pi_B^2} \, u = \sgn(\Pi_B)\,\Pi_B\,u$ and 
$\Pi_B\,u=\sgn(\Pi_B)\sqrt{\Pi_B^2} \, u$. 
\hfill$\Box$

\section{The homogeneous problem} \label{Ssubkato3}

In this section we prove the homogeneous subelliptic Kato problem
by an application of the results in Axelsson--Keith--McIntosh \cite{AKM} and Morris \cite{Morris1}
(see also Bandara \cite{Band1}).
Let $(\cx,d,\mu)$ be the metric measure space and let $(\Gamma,B_1,B_2)$
be the triple of operators in the Hilbert space
$\ch$ as in Subsection~\ref{Ssubkato2.2}.
We recall that for a scalar valued Lipschitz function $\eta$
we denote by $\Lipp\eta$ the pointwise Lipschitz constant.
If no confusion is possible, then we identify a measurable function
with the associated multiplication operator.

The hypothesis that are required are as follows.

\begin{enumerate}
\item[(H4)]
The metric space $(\cx,d)$ is complete and connected.
The measure $\mu$ is Borel-regular and doubling.
Moreover, there exists a number $N \in \Ni$ such that 
$\ch = L_2(\cx,\Ci^N;\mu)$.

\item[(H5)]
The operators $B_1$ and $B_2$ are multiplication operators 
by bounded matrix valued functions, denoted again by 
 $B_1,B_2 \in L_\infty(\cx,\cl(\Ci^N))$.

\item[(H6)]
There exists a constant $C > 0$ such that for every 
bounded Lipschitz function $\eta \colon \cx \to \Ri$ one has
\begin{enumerate}
\item 
the multiplication operator $\eta \, I$ leaves $D(\Gamma)$ invariant,
and,
\item
the commutator $[\Gamma, \eta \, I]$ is again a multiplication 
operator satisfying the bounds
\[
|([\Gamma, \eta \, I] u)(x)|
\leq C \, |(\Lipp \eta)(x)| \, |u(x)|
\]
for a.e.\ $x \in \cx$ and all $u \in D(\Gamma)$.
\end{enumerate} 

\item[(H7)]
For every open ball $B$ in $\cx$ one has
\[
\int_B \Gamma u \, d\mu = 0
\quad \mbox{and} \quad
\int_B \Gamma^* v \, d\mu = 0
\]
for all $u \in D(\Gamma)$ and $v \in D(\Gamma^*)$
with $\supp u \subset B$ and $\supp v \subset B$. 

\item[(H8)]
There exist $C_1,C_2 > 0$, $M \in \Ni$ and an operator 
$ Z  \colon D( Z ) \subset L_2(\cx,\Ci^N) \to L_2(\cx,\Ci^M)$, such that 
\begin{enumerate}
\item 
$D(\Pi) \cap \ran(\Pi) \subset D( Z )$,
\item (coercivity) \quad
$\| Z  u\| \leq C_2 \, \|\Pi u\|$  \\[5pt]
for all $u \in D(\Pi) \cap \ran(\Pi)$, and, 
\item (Poincar\'e estimate) \quad
$\displaystyle \int_B |u - \langle u \rangle_B|^2 \, d\mu
\leq C_1 \, r^2 \, \int_B | Z  u|^2 \, d\mu$  \\[5pt]
for all $x \in \cx$, $r > 0$ and $u \in D(\Pi) \cap \ran(\Pi)$, where
$B = B(x,r)$  and $\langle u \rangle_B:=\frac1{\mu(B)}\int_B u \, d\mu$. 
\end{enumerate}
\end{enumerate}

Hypothesis (H6) implies that $\Gamma$ behaves like a first order
differential operator.

The required quadratic estimates for the operator $\Pi_B$ 
now follow almost from Theorem~2.4 in Bandara \cite{Band1}.

\begin{thm} \label{tsubkato302}
Suppose $(\cx,d,\mu)$, $\ch$ and the triple $(\Gamma,B_1,B_2)$
satisfy Hypotheses {\rm (H1)--(H8)}.
Then the operator $\Pi_B$ satisfies the quadratic estimate
\[
\int_0^\infty \|Q^B_t u\|^2 \, \frac{dt}{t}
\simeq \|u\|^2
\]
for all  $u \in \overline{\ran(\Pi_B)} \subset \ch$. 
Hence for all $\mu \in (\omega,\frac{\pi}{2})$, the operator
$\Pi_B$ has a bounded $H^\infty(S_\mu^o)$-functional calculus, 
where $\omega$ is as in~{\rm (\ref{eSsubkato2;12})}.
\end{thm}
\proof\
The change to (H8) only affects Proposition~5.9 in Bandara \cite{Band1}.
This change forces the weighted Poincar\'e inequality
Proposition~5.8 in \cite{Band1} to become
\[
\int_{\cx} |u(x) - u_Q|^2 \left\langle \frac{d(x,Q)}{t} \right\rangle^{-M} \, d\mu(x)
	\lesssim \int_{\cx} |t \, ( Z  u)(x)|^2 \left\langle \frac{d(x,Q)}{t}\right\rangle^{p - M}
        \, d\mu(x)
\]
for all $u \in \ran(\Pi) \cap D(\Pi)$.
Consequently, in the proof of Proposition~5.9 in \cite{Band1}
we can invoke our coercivity Hypothesis (H8)
similar to the proof of Proposition~5.5 in \cite{AKM} to obtain the same conclusion.
The rest of the proof of Theorem~2.4 in \cite{Band1} remains unchanged.\hfill$\Box$

\vertspace

For the proof of Theorem~\ref{tsubkato102} we apply Theorem~\ref{tsubkato302}. 
Recall that in this case,  $G$ is the local direct product of a connected compact Lie group 
and a connected nilpotent Lie group.

\vertspace
\noindent
{\bf Proof of Theorem~\ref{tsubkato102}\ }\
Choose $\cx = G$, with Haar measure and distance being the subelliptic distance.
Let $\ch = L_2(G) \oplus L_2(G,\Ci^m) = L_2(G,\Ci^{1+m})$.
Define $D(\Gamma) = W_{1,2}'(G) \oplus L_2(G,\Ci^m)$ and define $\Gamma \colon D(\Gamma) \to \ch$
by 
\[
\Gamma = \left( \begin{array}{@{}cc@{}}
    0 & 0  \\[5pt]
    \nabla & 0 
                 \end{array} \right)
.  \]
Next let $B$ be the multiplication operator on $L_2(G,\Ci^m)$
by bounded matrix valued functions $(b_{kl})$.
Define $B_1,B_2 \colon \ch \to \ch$ by 
\[
B_1 = \left( \begin{array}{@{}cc@{}}
    b & 0  \\[5pt]
    0 & 0 
                  \end{array} \right)
\quad \mbox{and} \quad
B_2 = \left( \begin{array}{@{}cc@{}}
    0 & 0  \\[5pt]
    0 & B 
                  \end{array} \right)
.  \]
Then 
\[
\Gamma^* = \left( \begin{array}{@{}cc@{}}
    0 & -\divv  \\[5pt]
    0 & 0 
                  \end{array} \right)
\quad \mbox{and} \quad
\Pi = \left( \begin{array}{@{}cc@{}}
    0 & -\divv  \\[5pt]
    \nabla & 0 
                  \end{array} \right)
.  \]
Note that $\Pi$ is self-adjoint.
We next verify that Hypotheses (H1)--(H8) are valid.

\begin{enumerate}
\item[(H1)]
Since $C_c^\infty(G)$ and therefore also 
$W_{1,2}'(G)$ is dense in $L_2(G)$, and all the $A_k$ are closed operators,
this hypothesis is obvious.

\item[(H2)]
Clearly $B_1$ and $B_2$ are bounded. 
First note that 
\[
\RRe (B_1(f,0), (f,0))
 =  \int_G b \, |f|^2  
 \geq  \kappa_1 \, \|f\|_{L_2(G)}^2  
= \kappa_1 \, \|(f,0)\|_\ch^2
\] for all $f\in L^2(G)$, and hence when $f=-\divv w$ for all $w\in D(\divv)$.

Second, if $w \in W_{1,2}'(G)$, then
\[
\RRe (B_2(0,\nabla w),(0,\nabla w))
= \RRe \int_G \langle B \nabla w, \nabla w \rangle
\geq \kappa \, \int_G |\nabla w|^2
= \kappa \, \|(0,\nabla w)\|_\ch^2
,  \]
where we used the G\aa rding inequality (\ref{homgarding}).

\item[(H3)]
This trivially holds, since $B_1 \, B_2 = B_2 \, B_1 = 0$.

\item[(H4)]
The metric space $(G,d)$ is complete by Proposition~\ref{psubkato203}.\ref{psubkato203-1}.
By assumption, $G$ is connected.
The Haar measure is Borel-regular by Sections~11 and 15 in \cite{HR1}.
Since $G$ is the local direct product of a connected compact Lie group 
and a connected nilpotent Lie group it has polynomial growth.
Therefore the metric measure space $(G,d,\mu)$ has the doubling property by 
Proposition~\ref{psubkato203.5}.\ref{psubkato203.5-4}.

\item[(H5)]
This is obvious.

\item[(H6)]
Let $\eta$ be a bounded real valued Lipschitz function on $G$.
Then $\eta \in \bigcap_{k=1}^m D(A^{(\infty)}_k)$ by Proposition~\ref{psubkato210}.\ref{psubkato210-1}.
Let $u = (u_1,u_2) \in D(\Gamma)$.
Then $u_1 \in W_{1,2}'(G)$ and therefore $\eta \, u_1 \in W_{1,2}'(G)$.
So $\eta \, u \in D(\Gamma)$.
Moreover, for a.e.\ $x \in G$ one has
\begin{eqnarray*}
([\Gamma, \eta \, I] u)(x)
& = & (0,(\nabla (\eta \, u_1))(x)) - (0, (\eta \, \nabla u_1)(x))  \\
& = & (0, (A_1^{(\infty)} \eta)(x) \, u_1(x), \ldots, (A_m^{(\infty)} \eta)(x) \, u_m(x))
.  
\end{eqnarray*}
Hence $|([\Gamma, \eta \, I] u)(x)| \leq (\Lipp \eta)(x) \, |u(x)|$
by Proposition~\ref{psubkato210}.\ref{psubkato210-2}.

\item[(H7)]
Let $B$ be an open ball in $G$ and let $u = (u_1,u_2) \in D(\Gamma)$
with $\supp u \subset B$.
There exists a function $\chi \in C_c^\infty(B)$ such that $\chi(x) = 1$ for all 
$x \in \supp u$.
Then 
\begin{eqnarray*}
\int_B \Gamma u
= (0, \int_B \nabla u_1)
= (0, \int_B (\nabla u_1) \, \overline \chi)  
& = & (0, (A_1 u_1,\chi), \ldots, (A_m u_1, \chi))  \\
& = & - (0, (u_1,A_1 \chi), \ldots, (u_1,A_m \chi))  
= (0,0)
.
\end{eqnarray*}
Similarly, let $v = (v_1,v_2) \in D(\Gamma^*)$
with $\supp v \subset B$.
There exists a function $\chi \in C_c^\infty(B)$ such that $\chi(x) = 1$ for all 
$x \in \supp v$.
Then 
\begin{eqnarray*}
\int_B \Gamma^* v
= (- \int_B \divv v_2, 0)
= ( - \int_B \divv v_2 \, \overline \chi, 0)
& = & ( - (\divv v_2, \chi)_{L_2(G)}, 0)  \\
& = & ( ( v_2, \nabla \chi)_{L_2(G,\Ci^m)}, 0)
= (0,0)
.
\end{eqnarray*}

\item[(H8)]
Define the operator 
$ Z  \colon D( Z ) \to L_2(G,\Ci^m) \oplus L_2(G,\Ci^{m^2})) \simeq L_2(G,\Ci^{m + m^2})$
by $D( Z ) = W_{1,2}'(G) \oplus W_{1,2}'(G,\Ci^m) \subset \ch$ and 
\[
 Z (u_1,u_2) = (\nabla u_1,\widetilde \nabla u_2)
\]
where $\widetilde\nabla$ is defined in (\ref{tildenabla}). 
Let $(u_1,u_2) \in D(\Pi) \cap \ran(\Pi)$.
Then there exists an element $(v_1,v_2) \in D(\Pi)$ such that $(u_1,u_2) = \Pi(v_1,v_2)$.
This implies that $u_1,v_1 \in W_{1,2}'(G)$, $u_2,v_2 \in D(\divv)$ and, moreover, 
$(u_1,u_2) = (-\divv v_2, \nabla v_1)$.
Therefore $\nabla v_1 = u_2 \in D(\divv)$ and hence 
\[
v_1 \in D(\divv \nabla) 
= D(\Delta)
= \bigcap_{k,l=1}^m D(A_k \, A_l)
\]
by Proposition~\ref{psubkato202}.
So $u_2 = \nabla v_1 \in W_{1,2}'(G,\Ci^m)$.
This implies that $(u_1,u_2) \in D( Z )$.
We proved that $D(\Pi) \cap \ran(\Pi) \subset D( Z )$.
Moreover, if $C_2 \geq 1$ is as in Proposition~\ref{psubkato202} then 
\begin{eqnarray*}
\| Z (u_1,u_2)\|^2
= \|\nabla u_1\|^2 + \|\widetilde \nabla u_2\|^2  
& = & \|\nabla u_1\|^2 + \sum_{k,l=1}^m \|A_k \, A_l v_1\|^2  \\
& \leq & C_2 ( \|\nabla u_1\|^2 + \|\Delta v_1\|^2 )  \\
& = & C_2 ( \|\nabla u_1\|^2 + \|\divv u_2\|^2 )  \\
& = & C_2 \, \|\Pi(u_1,u_2)\|^2
.
\end{eqnarray*}
Finally, the Poincar\'e estimate follows from (P.1) in Saloff-Coste--Stroock \cite{SS} (page 118).
\end{enumerate}

Now it follows from Theorem~\ref{tsubkato302} that the operator $\Pi_B$ 
has a bounded $H^\infty$ functional calculus.
Since $\sgn \in \Hol^\infty(S_\mu^o)$ for all $\mu \in (0,\frac{\pi}{2})$, 
one has $D(\Pi_B) = D(\sqrt{\Pi_B^2})$ and 
$\|\Pi_B u\| \simeq \|\sqrt{\Pi_B^2} u\|$ for all $u \in D(\Pi_B)$
by Theorem~\ref{tsubkato211}.\ref{tsubkato205-2}.
Note that 
\[
\Pi_B = \left( \begin{array}{@{}cc@{}}
    0 & -b \divv B \\[5pt]
    \nabla & 0 
                  \end{array} \right)
\quad \mbox{and} \quad
\Pi_B^2 = \left( \begin{array}{@{}cc@{}}
    b H & 0 \\[5pt]
    0 & \widetilde H
                  \end{array} \right)
,  \]
where $\widetilde H u = - \nabla (b \divv (Bu))$ for all $u \in L_2(G,\Ci^m)$ 
with $(0,u) \in D(\Pi_B^2)$.
Then 
\[
\sqrt{\Pi_B^2} = \left( \begin{array}{@{}cc@{}}
    \sqrt{b H} & 0 \\[5pt]
    0 & \sqrt{\widetilde H}
                  \end{array} \right)
\]
and $\|\Pi_B (u_1,u_2)\|^2 = \|\nabla u_1\|^2 + \|b \divv (B u_2)\|^2$
for all $(u_1,u_2) \in D(\Pi_B) = W_{1,2}'(G) \oplus D(\divv \circ B)$.
Restricting to the scalar valued functions gives
$D(\sqrt{b H}) = W_{1,2}'(G)$ and $\|\sqrt{b H} \, u\| \simeq \|\nabla u\|$.
This completes the proof of Theorem~\ref{tsubkato102}.\hfill$\Box$

\vertspace

We conclude this section with a proof of the homogeneous 
estimates of Theorem~\ref{tsubkato303} for 
homogeneous strongly elliptic operators on connected Lie groups 
with polynomial growth.

\vertspace
\noindent
{\bf Proof of Theorem~\ref{tsubkato303}\ }\
We use the structure theory for Lie groups with polynomial growth as developed 
in Dungey--ter Elst--Robinson \cite{DER4} and summarised on pages 125--126.
There exists another group multiplication $*$ on $G$ such that the manifold
$G$ with multiplication $*$ is a Lie group, denoted by $G_N$, which is 
the local direct product of a connected compact Lie group 
and a connected nilpotent Lie group.
Let $\gotg_N$ be the Lie algebra of $G_N$.
Then $\gotg_N = \gotg$ as vector spaces.
The Haar measure $\mu$ on $G$ is again a Haar measure on $G_N$.
Moreover, $W_{1,2}(G,a) = W_{1,2}(G_N,a)$.
There exist a Lie group homomorphism $\overline{\cs} \colon G_N \to \Aut(\gotg_N)$,
an inner product $\langle \cdot , \cdot \rangle$ on $\gotg_N$ and an
orthonormal basis $b_1,\ldots,b_m$ of $\gotg_N$ such that 
$\overline{\cs}(x)$ is orthogonal on $\gotg_N$ for all $x \in G$ and 
\[
(A_k u)(x)
= \sum_{l=1}^m \langle \overline{\cs}(x) \, b_l, a_k \rangle \, (B^{(N)}_l u)(x)
\]
for a.e.\ $x \in G$, $k \in \{ 1,\ldots,m \} $ and $u \in W_{1,2}(G,a)$,
where $B^{(N)}_l$ is the infinitesimal generator of $G_N$ in the direction $b_l$.
Note that $W_{1,2}(G,a) = W_{1,2}(G_N,a) = W_{1,2}(G_N,b)$
since both $a_1,\ldots,a_m$ and $b_1,\ldots,b_m$ are vector space bases
in $\gotg$. 
Moreover, there exists a constant $C \geq 1$ such that 
\begin{equation}
\frac{1}{C} \sum_{k=1}^m \|A_k u\|^2 
\leq \sum_{k=1}^m \|B^{(N)}_k u\|^2 
\leq C \sum_{k=1}^m \|A_k u\|^2 
\label{etsubkato303;1}
\end{equation}
for all $u \in W_{1,2}(G,a)$.

Now
\[
H = - \sum_{k',l' = 1}^m B^{(N)}_{k'} \, \tilde c_{k' \, l'} \, B^{(N)}_{l'}
\]
as a divergence form operator, where
\[
\tilde c_{k' \, l'}(x)
= \sum_{k,l=1}^m \langle \overline{\cs}(x) \, b_{k'} , a_k \rangle  \, 
     c_{kl}(x) \, \langle \overline{\cs}(x) \, b_{l'}, a_l \rangle
.  \]
Using (\ref{etsubkato303;1}) it follows that $H$ satisfies the 
assumptions of Theorem~\ref{tsubkato102} 
with respect to the coefficients $\tilde c_{kl}$, the group $G_N$ 
and the basis $b_1,\ldots,b_m$.
Therefore $\|\sqrt{b H} \, u\| \simeq \sum_{k=1}^m \|B^{(N)}_k u\|$ 
by Theorem~\ref{tsubkato102}.
A further use of the equivalence (\ref{etsubkato303;1}) of the seminorms,
completes the proof of the theorem.\hfill$\Box$

\section{The inhomogeneous problem} \label{Ssubkato4}

To solve the inhomogeneous problem, we apply the results on
quadratic estimates in the framework of Morris \cite{Morris1},
with appropriate modifications, in particular to Hypotheses (H6)
and (H8) of \cite{Morris1}.

For us, this means we continue to use Hypotheses (H1)--(H3) from 
Section~\ref{Ssubkato2}, and Hypotheses (H5), (H6) from Section~\ref{Ssubkato3}, 
together with the following ones. \

\begin{enumerate}

\item[(H4i)]
	The metric measure space $(\cx,d,\mu)$ is 
	a complete, connected, locally exponentially doubling 
	metric measure space with
	a Borel-regular measure $\mu$.
Moreover, there exists a number $N \in \Ni$ such that 
$\ch = L_2(\cx,\Ci^N;\mu)$.

\item[(H7i)]
There exists a constant $c > 0$ such that for every open ball $B$ in $\cx$
with radius at most~$1$ one has
\[
\Big| \int_B \Gamma u \, d\mu \Big| \leq c \, \mu(B)^{\frac{1}{2}} \, \|u\|
\quad \mbox{and} \quad
\Big| \int_B \Gamma^* v \, d\mu \Big| \leq c \, \mu(B)^{\frac{1}{2}} \, \|v\|
\]
for all $u \in D(\Gamma)$ and $v \in D(\Gamma^*)$
with $\supp u \subset B$ and $\supp v \subset B$. 

\item[(H8i)]
There exist $C_1,C_2 > 0$, $M \in \Ni$ and an operator 
$Z \colon D(Z) \to L_2(\cx,\Ci^M)$ with 
$D(Z) \subset L_2(\cx,\Ci^N) = \ch$, such that 
\begin{enumerate}
\item 
$D(\Pi) \cap \ran(\Pi) \subset D(Z)$,
\item (coercivity) \quad
$\|Z u\| + \|u\| \leq C_2 \, \|\Pi u\|$  \\[5pt]
for all $u \in D(\Pi) \cap \ran(\Pi)$, and,
\item (Poincar\'e estimate) \quad
$\displaystyle \int_B |u - \langle u \rangle_B|^2 \, d\mu
\leq C_1 \, r^2 \, \int_B (|Z u|^2 + |u|^2) \, d\mu$  \\[5pt]
for all $x \in \cx$, $r \in (0,\infty)$ and $u \in D(\Pi) \cap \ran(\Pi)$, where
$B = B(x,r)$.
\end{enumerate}
\end{enumerate} 

\begin{remarkn}\label{poincare} 
For any $r_0 > 0$, the Poincar\'e estimate in (H8i) is valid uniformly  
for all $r \geq r_0$. 
\end{remarkn} 

We are now able to formulate the theorem on quadratic estimates
for the inhomogeneous operators.

\begin{thm} \label{tsubkato401}
Suppose $(\cx,d,\mu)$, $\ch$ and the triple $(\Gamma,B_1,B_2)$
satisfy Hypotheses {\rm (H1), (H2), (H3), (H4i), (H5), (H6), (H7i)} 
and {\rm (H8i)}.
Then the operator $\Pi_B$ satisfies the quadratic estimate
\[
\int_0^\infty \|Q^B_t u\|^2 \, \frac{dt}{t}
\simeq \|u\|^2
\]
for all $u \in \overline{\ran(\Pi_B)} \subset L_2(\cx,\Ci^N)$.
Hence $\Pi_B$ has a bounded $H^\infty$-functional calculus.
\end{thm}
\proof\
First, we note that for all non-empty subsets $E, F$ on any metric space $X$
satisfying $d(E,F) > 0$
one can find a Lipschitz function
$\eta \colon X \to [0,1]$ such that $\eta = 1$ on~$E$,
$\eta = 0$ on~$F$, and ${\mathop{\rm\mathbf{Lip}}}\ \eta \leq 1/d(E,F)$,
where ${\mathop{\rm\mathbf{Lip}}}\ \eta$ is the Lipschitz constant of $\eta$.
Also observe that all the smooth cutoff functions used
in Morris \cite{Morris1} can be replaced by Lipschitz equivalents, 
in particular allowing us to obtain off-diagonal estimates in the present case.

Next, our alteration to (H6) and (H8) allows us to 
dispense with the use of the Sobolev spaces with respect to the
Levi--Civita connection in \cite{Morris1} and simply consider
$\ran(\Pi) \cap D(\Pi)$.
Explicitly, the weighted Poincar\'e inequality in Lemma~5.7
of \cite{Morris1} is altered to read $\|Z u\|$ instead of
$\|\nabla u\|$, and by coercivity, Proposition~5.8 in \cite{Morris1} holds.
Then Proposition~5.2 in \cite{Morris1} holds, 
since we still have $\|u\| \lesssim \|\Pi u\|$.
Finally, we observe that the measure merely needs to be locally exponentially 
doubling and Borel-regular.\hfill$\Box$

\vertspace

Now we are able to prove Theorem~\ref{tsubkato101}, using ideas from~\cite{Morris1},
which have their roots in \cite{AKM2}. 

\vertspace
\noindent
{\bf Proof of Theorem~\ref{tsubkato101}\ }\
We apply Theorem~\ref{tsubkato401} with $\cx = G$ and 
$\ch = L_2(G) \oplus \big( L_2(G) \oplus L_2(G,\Ci^m) \big) = L_2(G,\Ci^{2+m})$.
Define $D(\Gamma) = W_{1,2}'(G) \oplus L_2(G) \oplus L_2(G,\Ci^m)$ 
and define $\Gamma \colon D(\Gamma) \to \ch$
by 
\[
\Gamma = \left( \begin{array}{@{}ccc@{}}
    0 & 0 & 0  \\[5pt]
    I & 0 & 0  \\[5pt]
    \nabla & 0 & 0
                 \end{array} \right)
.  \]
Let $B$ be the multiplication operator on $L_2(G) \oplus L_2(G,\Ci^m) = L_2(G,\Ci^{1+m})$
by bounded matrix valued functions
\[
\left( \begin{array}{@{}cccc@{}}
   b_0 & b_1' & \cdots & b_m'  \\[5pt]
   b_1 & & & \\[5pt]
   \vdots & & (b_{kl}) &  \\[5pt]
   b_m & & & 
\end{array} \right)
.  \]
Next define $B_1,B_2 \colon \ch \to \ch$ by 
\[
B_1 = \left( \begin{array}{@{}cc@{}}
    b & 0  \\[5pt]
    0 & 0 
                  \end{array} \right)
\quad \mbox{and} \quad
B_2 = \left( \begin{array}{@{}cc@{}}
    0 & 0  \\[5pt]
    0 & B 
                  \end{array} \right)
.  \]
Then 
\[
\Gamma^* = \left( \begin{array}{@{}ccc@{}}
    0 & I & -\divv  \\[5pt]
    0 & 0 & 0  \\[5pt]
    0 & 0 & 0 
                  \end{array} \right)
\quad \mbox{and} \quad
\Pi = \left( \begin{array}{@{}ccc@{}}
    0 & I & -\divv  \\[5pt]
    I & 0 & 0  \\[5pt]
    \nabla & 0 & 0 
                  \end{array} \right)
.  \]
Note again that $\Pi$ is self-adjoint.

The proof that Hypotheses (H1), (H3), (H5) and (H6) are valid
is similar to the proof of Theorem~\ref{tsubkato102} in Section~\ref{Ssubkato3}.
Also (H2) follows similarly from the G\aa rding inequality (\ref{eSsubkato1;20}).
It remains to verify Hypotheses (H4i), (H7i) and (H8i).

\begin{enumerate}
\parindent=20pt 
\parskip=0pt
\item[(H4i)]
The only difference with (H4) is the locally exponentially doubling
property, which follows from Proposition~\ref{psubkato203.5}.\ref{psubkato203.5-1}.

\item[(H7i)]
Let $B$ be an open ball in $G$ and let $u = (u_1,u_2,u_3) \in D(\Gamma)$
with $\supp u \subset B$.
Then $\int_B \nabla u_1 = 0$ by the argument in the proof of (H7) in Section~\ref{Ssubkato3}.
Therefore
\[
\Big| \int_B \Gamma u \Big|^2
= \Big| \int_B u_1 \Big|^2 + \Big| \int_B \nabla u_1 \Big|^2
= \Big| \int_B u_1 \Big|^2
\leq \mu(B) \, \|u_1\|^2
\leq \mu(B) \, \|u\|^2
.  \]
Similarly, if $v = (v_1,v_2,v_3) \in D(\Gamma^*)$ with $\supp v \subset B$ then 
$\int_B \divv v_3 = 0$ and 
\[
\Big| \int_B \Gamma^* v \Big|
= \Big| \int_B v_2 - \divv v_3 \Big|
= \Big| \int_B v_2 \Big|
\leq \mu(B)^{1/2} \, \|v\|
\]
as required.

\item[(H8i)]
Define the operator 
$Z \colon D(Z) \to 
L_2(G,\Ci^m) \oplus L_2(G,\Ci^m) \oplus L_2(G,\Ci^{m^2}) = L_2(G,\Ci^{2m + m^2})$
by $D(Z) = W_{1,2}'(G) \oplus W_{1,2}'(G) \oplus W_{1,2}'(G,\Ci^m) \subset \ch$ and 
\[
Z(u_1,u_2,u_3) = (\nabla u_1,\nabla u_2,\widetilde \nabla u_3)
.  \]
Let $u=(u_1,u_2,u_3) \in D(\Pi) \cap \ran(\Pi)$. 
Then there exists an element $(v_1,v_2,v_3) \in D(\Pi)$ such that $(u_1,u_2,u_3) = \Pi(v_1,v_2,v_3) =(v_2-\divv v_3, v_1, \nabla v_1)$
with in particular, $u_1,u_2=v_1 \in W_{1,2}'(G)$ and $u_3,v_3 \in D(\divv)$.
Therefore $\nabla v_1 = u_3 \in D(\divv)$ and hence 
\[
v_1 \in D(\divv \nabla) 
= D(\Delta)
= \bigcap_{k,l=1}^m D(A_k \, A_l)
\]
by Proposition~\ref{psubkato202}. 
So  $u_3 = \nabla v_1 \in W_{1,2}'(G,\Ci^m)$.
This implies that  $u=(u_1,u_2,u_3) \in D(Z)$.  
We proved that $D(\Pi) \cap \ran(\Pi) \subset D(Z)$,  
and it remains for us to obtain the bound $\|Z u \| + \|u\|\leq C_2\|\Pi u\|$. 

Note that 
\[
\|Z u\|^2
= \|\nabla u_1\|^2 + \|\nabla u_2\|^2 + \|\widetilde \nabla u_3\|^2  
\]
and
\[
\|\Pi u\|^2
= \|u_2 - \divv u_3\|^2 + \|u_1\|^2 + \|\nabla u_1\|^2
= \|(I - \Delta) v_1\|^2 + \|u_1\|^2 + \|\nabla u_1\|^2
.  \] 
Clearly 
\[
\|\nabla u_2\|^2
= \|\nabla v_1\|^2
= (\Delta v_1, v_1)
\leq \|\Delta v_1\| \, \|v_1\|
\leq \|(I - \Delta) v_1\|^2
\]
and if $C_1 \geq 1$  as in Proposition~\ref{psubkato202},  then 
\[
\|\widetilde \nabla u_3\|^2  
= \sum_{k,l=1}^m \|A_k \, A_l v_1\|^2
\leq C_1 (\|\Delta v_1\|^2 + \|v_1\|^2)  
\leq 2 C_1 \|(I - \Delta) v_1\|^2
.  \]
So 
\[
\|Z(u_1,u_2,u_3)\|^2
\leq 3 C_1 \|\Pi(u_1,u_2,u_3)\|^2
.  \]
Next we show that $\|(u_1,u_2,u_3)\|^2 \leq (4 C_1 \, m^2 + 4) \|\Pi(u_1,u_2,u_3)\|^2$.
Trivially, $\|u_1\|^2 \leq \|\Pi(u_1,u_2,u_3)\|^2$.
Moreover, 
\[
\|\divv u_3\|^2 \leq m^2 \|\widetilde \nabla u_3\|^2
\leq 2 C_1 \, m^2 \|(I - \Delta) v_1\|^2
\]
and 
\[
\|u_2\|^2 
\leq 2 \|u_2 - \divv u_3\|^2 + 2 \|\divv u_3\|^2
\leq (4 C_1 \, m^2 + 2) \|(I - \Delta) v_1\|^2
.  \]
Also $\|u_3\|^2 = \|\nabla v_1\|^2 \leq \|(I - \Delta) v_1\|^2$.

We conclude that $\|Z u \| + \|u\|\leq C_2\|\Pi u\|$ 
(for a suitable constant $C_2$) 
as required. 

Finally we prove the Poincar\'e inequality.
Jerison \cite{Jer} proved that there exist $c,r_0 > 0$ such that 
\begin{equation} 
\int_{B(r)} |f - \langle f \rangle_{B(r)}|^2
\leq c \, r^2 \int_{B(r)} |\nabla f|^2
\label{etsubkato101;1}
\end{equation}
for all $r \in (0,r_0]$ and $f \in C^\infty(\overline{B(r)})$.
Since $C^\infty(G) \cap W_{1,2}'(G)$ is dense in $W_{1,2}'(G)$ by 
Lemma~2.4 in ter Elst--Robinson \cite{ER5}, it follows that (\ref{etsubkato101;1}) is 
valid for all $r \in (0,r_0]$ and $f \in W_{1,2}'(G)$.
Let $R$ denote the right regular representation in $G$
and let $\delta$ be again the modular function on $G$.
Then 
\begin{eqnarray*}
\int_{B(x,r)} |f - \langle f \rangle_{B(x,r)}|^2
& = & \delta(x) \, \int_{B(r)} |R(x) f - \langle R(x)f \rangle_{B(r)}|^2  \\
& \leq & c \, \delta(x) \, r^2 \int_{B(r)} |\nabla \, R(x) f|^2
= c \, r^2 \int_{B(x,r)} |\nabla f|^2
\end{eqnarray*}
for all $x \in G$, $r \in (0,r_0]$ and $f \in W_{1,2}'(G)$.
The Poincar\'e estimate follows by  Remark~\ref{poincare}. 
\end{enumerate}

Now one can complete the proof of Theorem~\ref{tsubkato101}
similarly as in Section~\ref{Ssubkato3} by an application
of Theorem~\ref{tsubkato401}.\hfill$\Box$

\section{Further results} \label{Ssubkato5}

We have chosen to take the (left) Haar measure $\mu$ on $L_2(G;\mu)$,
and the infinitesimal generators are with respect to the left regular
representation in $L_2(G;\mu)$.
Another option would be to choose the right Haar measure $\nu$ on $G$
and consider the left regular representation in $L_2(G;\nu)$.
Then the solution to the Kato problem has the following formulation.

\begin{thm} \label{tsubkato402}
Let $a_1,\ldots,a_m$ be an algebraic basis for the Lie algebra $\gotg$ of a 
connected Lie group~$G$.
For all $k \in \{ 1,\ldots,m \} $ let $A^{(R)}_k$ be the infinitesimal 
generator of the one-parameter group $t \mapsto L^{(R)}(\exp t a_k)$,
where $L^{(R)}$ denotes the left regular representation in $L_2(G;\nu)$.
For all $k,l \in \{ 1,\ldots,m \} $ let $b_{kl}, b_k,b_k',b_0 \in L_\infty(G)$.
Assume there exist $\kappa,c_1 > 0$ such that 
\[
\RRe \sum_{k,l=1}^m (b_{kl} \, A^{(R)}_k u, A^{(R)}_l u) 
\geq \kappa \sum_{k=1}^m \|A^{(R)}_k u\|^2 - c_1 \, \|u\|^2
\]
for all $u \in \bigcap_{k=1}^m D(A^{(R)}_k)$.
Consider the divergence form operator
\[
H = \sum_{k,l=1}^m (A^{(R)}_k)^* \, b_{kl} \, A^{(R)}_l 
   + \sum_{k=1}^m b_k \, A^{(R)}_k
   + \sum_{k=1}^m (A^{(R)}_k)^* \, b_k'
   + b_0 \, I
\]
in $L_2(G;\nu)$, where the norm and inner product are in $L_2(G;\nu)$.
Suppose $\RRe b_0$ is large enough such that $- H$
generates a bounded semigroup on $L_2(G;\nu)$.
Let $b \in L_\infty(G)$ and suppose there exists a constant $\kappa_1 > 0$ such that 
$\RRe b \geq \kappa_1$ a.e.
Then 
\[
D(\sqrt{b H}) = \bigcap_{k=1}^m D(A^{(R)}_k)
\]
with equivalent norms.
\end{thm}
\proof\
The proof is almost the same as the proof of Theorem~\ref{tsubkato101},
so we indicate the differences. 
We replace all $A_k$ by $A^{(R)}_k$.
Note that for all $k \in \{ 1,\ldots,m \} $ there exists a constant $\beta_k \in \Ri$
such that $(A^{(R)}_k)^* = - A^{(R)}_k + \beta_k \, I$, where 
the adjoint is in $L_2(G;\nu)$.
We take the same subelliptic distance on $G$.
 There exists a constant
$c > 0$ such that 
$\nu = c \, \delta^{-1} \mu$, where $\delta$ is the modular function.
(See \cite{HR1} Theorem~15.15.)
Moreover, since $\delta$ is a continuous homomorphism,
there exist $M,\omega > 0$ such that 
$\delta(x) \leq M \, e^{\omega d(x,e)}$
for all $x \in G$.
Hence by Proposition~\ref{psubkato203} there are $c,C,\lambda > 0$ and $D' \in \Ni$ 
such that $c \, r^{D'} \leq \nu(B(r)) \leq C \, r^{D'}$ for all $r \in (0,1]$
and $\nu(B(r)) \leq C \, e^{\lambda r}$ for all $r \geq 1$.
(Actually, the natural number $D'$ is the same as in 
Proposition~\ref{psubkato203}.\ref{psubkato203-2}.)
Then Hypothesis (H4i) follows.
Next consider Hypothesis (H7i).
Let $B$ be an open ball in $G$, let $u = (u_1,u_2,u_3) \in D(\Gamma)$
and suppose that $\supp u \subset B$.
There exists a function $\chi \in C_c^\infty(B)$ such that $\chi(x) = 1$ 
for all $x \in \supp u$.
Let $k \in \{ 1,\ldots,m \} $.
Then 
\[
\int_B A^{(R)}_k u_1 \, d\nu
= (u_1, (A^{(R)}_k)^* \, \overline \chi)_{L_2(G;\nu)}
= (u_1, ( -A^{(R)}_k + \beta_k \, I) \, \overline \chi)_{L_2(G;\nu)}
= \beta_k \int_B u_1 \, d\nu
.  \]
So $|\int_B \nabla u_1 \, d\nu| \leq \sqrt{\beta_1^2 + \ldots + \beta_m^2} \, \nu(B)^{1/2} \, \|u\|$.
The rest of the proof of Hypothesis (H7i) is similar.

All other hypothesis have the same proof as before.\hfill$\Box$

\vertspace

One can also consider the infinitesimal generators with respect to the 
right regular representation on $L_2(G;\mu)$ or $L_2(G;\nu)$.
Then the inhomogeneous Kato problem has again a solution.
This follows from Theorems~\ref{tsubkato101} and \ref{tsubkato402}
by using the inversion $x \mapsto x^{-1}$ on $G$.
We leave the formulation of the two theorems to the reader.

\section{Stability} \label{Ssubkato6}

Finally we consider stability under holomorphic perturbation.

Let $U \subset \Ci$ be an open set, $\omega \in [0,\frac{\pi}{2})$  
and for all $\zeta \in U$ let 
$T(\zeta)$ be an $\omega$-bisectorial operator in $\ch$ with domain $D(T(\zeta))$.
Let $\mu \in (\omega,\frac{\pi}{2})$.
We say that $T$ has a {\bf uniformly bounded holomorphic $H^\infty(S_\mu^o)$-functional calculus}
if there exists a constant $C > 0$ such that 
$\|\psi(T(\zeta))\| \leq C \, \|\psi\|_\infty$ uniformly for all 
$\psi \in \Psi(S_\mu^o)$ and $\zeta \in U$.

\begin{thm} \label{tsubkato206}
Let $U \subset \Ci$ be an open set,
$\ch$ a Hilbert space and $(X,d,\mu)$ a metric measure space.
Let $B_1,B_2 \colon U \to \cl(\ch)$
be bounded holomorphic functions.
Suppose that the triple $(\Gamma,B_1(\zeta),B_2(\zeta))$ satisfies 
{\rm (H1)--(H8)} uniformly for all $\zeta \in U$, with constants $\kappa_1$ and 
$\kappa_2$.
Let 
\[
\omega = \sup_{\zeta \in U}
\frac{1}{2} \Big( \arccos \frac{\kappa_1}{\|B_1(\zeta)\|} + \arccos \frac{\kappa_2}{\|B_2(\zeta)\|} \Big)
< \frac{\pi}{2}
.  \]
Let $\mu \in (\omega,\frac{\pi}{2})$.
Then one has the following.
\begin{tabel}
\item \label{tsubkato206-1}
The operator $\Pi_{B(\zeta)}$ is $\omega$-bisectorial operator in $\ch$
uniformly for all $\zeta \in U$.
\item \label{tsubkato206-2}
The family $\zeta \mapsto \Pi_{B(\zeta)}$ has a uniformly bounded holomorphic 
$H^\infty(S_\mu^o)$-functional calculus.
\item \label{tsubkato206-3}
For all $f \in \Hol^\infty(S_\mu^o)$ 
the map $\zeta \mapsto f(\Pi_{B(\zeta)})$ is holomorphic.
\end{tabel}
\end{thm}
\proof\
Statement~\ref{tsubkato206-1} follows from Proposition~2.5 in Axelsson--Keith--McIntosh \cite{AKM}
and Statement~\ref{tsubkato206-2} from Theorem~\ref{tsubkato302}.
Statement~\ref{tsubkato206-3} follows as in the proof of Theorem~6.4 in \cite{AKM}.\hfill$\Box$

\vertspace

We conclude the paper by noting the following stability result.

\begin{thm} \label{tsubkato501}
Let $G$ be the local direct product of a connected compact Lie group 
and a connected nilpotent Lie group.
Let 
\[
H = - \sum_{k,l=1}^m A_k \, b_{kl} \, A_l 
\]
be a homogeneous divergence form operator with bounded measurable coefficients
satisfying the subellipticity condition {\rm (\ref{homgarding})} with constant $\kappa_1$.
Let $b \in L_\infty(G)$ and suppose there exists a constant $\kappa_2 > 0$ such that 
$\RRe b \geq \kappa_2$ a.e.
Let $\eta_1 \in (0,\kappa_1)$ and $\eta_2 \in (0,\kappa_2)$.
Then there exists a constant $C > 0$ such that the following is valid.
For all $k,l \in \{ 1,\ldots,m \} $ let $\tilde b_{kl} \in L_\infty(G)$
and suppose that
\[
\widetilde M 
= \sup_{x \in G} \|(\tilde b_{kl}(x))\|_{\Ci^{m \times m}}
\leq \eta_1
.  \] 
Further, let $\tilde b \in L_\infty(G)$ with $\|\tilde b\|_\infty \leq \eta_2$.
Let 
\[
\widetilde H = - \sum_{k,l=1}^m A_k \, (b_{kl} + \tilde b_{kl}) \, A_l 
.  \]
Then 
\[
\|\sqrt{(b + \tilde b) \widetilde H} u - \sqrt{b H} u\| 
\leq C (\widetilde M + \|\tilde b\|_\infty) \|\nabla u\|
\]
for all $u \in W_{1,2}'(G)$.
\end{thm}
\proof\
This follows as in the proof of Theorem~6.5 in \cite{AKM}, using
Theorems~\ref{tsubkato302} and Theorem~\ref{tsubkato206} with 
$B(\zeta)=B+\zeta \tilde B$ for all $\zeta\in U$, where $U$ is an appropriate open 
set with $[0,1]\subset U\subset \C$.
See also the proof of Theorem~7.2 in \cite{BMC}.\hfill$\Box$

\vertspace

There are similar stability results for the inhomogeneous problems
as in Theorem~\ref{tsubkato101}, or with the right Haar measure or 
right translations. 
We leave the formulation to the reader.

\subsection*{Acknowledgements}
The authors appreciate the support of the Centre for Mathematics and its Applications 
at the Australian National University, Canberra, where this project was begun, 
the Department of Mathematics at the University of Auckland, and the Department of
Mathematics, University of Missouri, where it was continued.
The first author was supported through an Australian Postgraduate Award,
through the Mathematical Sciences Institute, Australian National University, Canberra,
and an Australian-American Fulbright Scholarship.
The second and third authors were supported by the Marsden Fund Council from N.Z. Government funding, 
administered by the Royal Society of New Zealand. 
The third author gratefully acknowledges support from the 
Australian Government through the Australian Research Council.

\newcommand{\etalchar}[1]{$^{#1}$}

\small 

\noindent
{\sc Lashi Bandara,
Centre for Mathematics and its Applications,
Mathematical Sciences Institute,
Australian National University,
Canberra, ACT 0200,
Australia}  \\
{\em E-mail address}\/: {\bf Lashi.Bandara@anu.edu.au}

\mbox{}

\noindent
{\sc A.F.M. ter Elst,
Department of Mathematics,
University of Auckland,
Private bag 92019,
Auckland 1142,
New Zealand}  \\
{\em E-mail address}\/: {\bf terelst@math.auckland.ac.nz}

\mbox{}

\noindent
{\sc Alan McIntosh,
Centre for Mathematics and its Applications,
Mathematical Sciences Institute,
Australian National University,
Canberra, ACT 0200,
Australia}  \\
{\em E-mail address}\/: {\bf Alan.McIntosh@anu.edu.au}

\end{document}